\documentclass[11 pt]{article}
\usepackage[T1]{fontenc}
\usepackage{stackengine, IEEEtrantools}

\usepackage{cite, authblk}
\usepackage{amsmath,amssymb,amsfonts}
\usepackage{algorithmic}
\usepackage{graphicx}
\usepackage{textcomp}
\usepackage{xcolor}
\usepackage{algorithm,algorithmic}
\usepackage{mathtools, stackengine, mathrsfs, setspace, bm, multirow, tikz, pgfplots, adjustbox, calc, url}
\usepackage{xparse}
\pgfplotsset{width=8.1cm, legend style={font=\small}}
\newsavebox{\fminipagebox}
\NewDocumentEnvironment{fminipage}{m O{\fboxsep}}
 {\par\kern#2\noindent\begin{lrbox}{\fminipagebox}
  \begin{minipage}{#1}\ignorespaces}
 {\end{minipage}\end{lrbox}%
  \makebox[#1]{%
    \kern\dimexpr-\fboxsep-\fboxrule\relax
    \fbox{\usebox{\fminipagebox}}%
    \kern\dimexpr-\fboxsep-\fboxrule\relax
  }\par\kern#2
 }
\begingroup\edef\x{\endgroup
  \mathchardef\mathdollar=\the\numexpr"7000+\the\mathdollar\relax
}\x
\DeclareMathAlphabet{\mathit}{T1}{cmr}{m}{it}

\makeatletter
\let\MYcaption\@makecaption
\makeatother
\usepackage{caption}
\DeclareCaptionLabelSeparator{colonquad}{:\quad}
\DeclareCaptionStyle{fig-top}%
	[justification=raggedright,indention=2cm, font=footnotesize, labelsep=colonquad]{}
\DeclareCaptionStyle{fig}%
	[justification=raggedright,indention=2cm, font=footnotesize, labelsep=colonquad]{}	
\makeatletter
\let\@makecaption\MYcaption
\makeatother

\pgfplotsset{
    every first x axis line/.style={},
    every first y axis line/.style={},
    every first z axis line/.style={},
    every second x axis line/.style={},
    every second y axis line/.style={},
    every second z axis line/.style={},
    first x axis line style/.style={/pgfplots/every first x axis line/.append style={#1}},
    first y axis line style/.style={/pgfplots/every first y axis line/.append style={#1}},
    first z axis line style/.style={/pgfplots/every first z axis line/.append style={#1}},
    second x axis line style/.style={/pgfplots/every second x axis line/.append style={#1}},
    second y axis line style/.style={/pgfplots/every second y axis line/.append style={#1}},
    second z axis line style/.style={/pgfplots/every second z axis line/.append style={#1}}
}

\makeatletter
\def\pgfplots@drawaxis@outerlines@separate@onorientedsurf#1#2{%
    \if2\csname pgfplots@#1axislinesnum\endcsname
    \else
    \scope[/pgfplots/every outer #1 axis line,
        #1discont,decoration={pre length=\csname #1disstart\endcsname, post length=\csname #1disend\endcsname}]
        \pgfplots@ifaxisline@B@onorientedsurf@should@be@drawn{0}{%
            \draw [/pgfplots/every first #1 axis line] decorate {
                \pgfextra
                \pgfplotspointonorientedsurfaceabsetupfor{#2}{#1}{\pgfplotspointonorientedsurfaceN}%
                \pgfplots@drawgridlines@onorientedsurf@fromto{\csname pgfplots@#2min\endcsname}%
                \endpgfextra 
                };
        }{}%
        \pgfplots@ifaxisline@B@onorientedsurf@should@be@drawn{1}{%
            \draw [/pgfplots/every second #1 axis line] decorate {
                \pgfextra
                \pgfplotspointonorientedsurfaceabsetupfor{#2}{#1}{\pgfplotspointonorientedsurfaceN}%
                \pgfplots@drawgridlines@onorientedsurf@fromto{\csname pgfplots@#2max\endcsname}%
                \endpgfextra 
                };
        }{}%
    \endscope
    \fi
}%
\makeatother
\newcommand\blfootnote[1]{%
  \begingroup
  \renewcommand\thefootnote{}\footnote{#1}%
  \addtocounter{footnote}{-1}%
  \endgroup
}
\def\BibTeX{{\rm B\kern-.05em{\sc i\kern-.025em b}\kern-.08em
    T\kern-.1667em\lower.7ex\hbox{E}\kern-.125emX}}

\allowdisplaybreaks

\title{Kullback-Leibler Divergence-Based \\Distributionally Robust Unit Commitment\\ Under Net Load Uncertainty
}
\author{Ogun Yurdakul}
\author{Fikret Sivrikaya}
\author{Sahin Albayrak}
\affil{Department of Electrical Engineering and Computer Science\\
Technical University of Berlin, Berlin, Germany\\
Email: \{yurdakul, fikret.sivrikaya, sahin.albayrak\}@tu-berlin.de}
\date{}                     
\begin{document}

\maketitle

\begin{abstract}
The deepening penetration of renewable resources into power systems entails great difficulties that have not been surmounted satisfactorily. An issue that merits special attention is the short-term planning of power systems under net load uncertainty. To this end, we work out a distributionally robust unit commitment methodology that expressly assesses the uncertainty associated with net load. The principal strength of the proposed methodology lies in its ability to represent the probabilistic nature of net load without having to set forth its probability distribution. This strength is brought about by the notion of ambiguity set, for the construction of which the Kullback-Leibler divergence is employed in this paper. We demonstrate the effectiveness of the proposed methodology on real-world data using representative studies. The sensitivity analyses performed provide quantitative answers to a broad array of \textit{what if} questions on the influence of divergence tolerance and dataset size on optimal solutions.\\\\
\end{abstract}
\blfootnote{This work was supported in part by the Research Council of Norway under the ``LUCS'' project, and by the German Federal Ministry for Economic Affairs and Energy under Grant 03EI6004B.}
\section{Introduction}\label{1}
The growing concerns for the environment bring about the enactment of technical, regulatory, and economic measures that aim at the deeper penetration of renewable resources (\textit{RR}s). These measures have begotten the use of the legacy grid in a way different from that for which it was designed, thus posing unprecedented challenges for grid operators (\textit{GO}s).\par
One such challenge is the short-term planning of power systems vis-à-vis uncertain net load, that is, system load less \textit{RR} generation. Short-term planning of power systems involves the unit commitment (\textit{UC}) and economic dispatch decisions. The \textit{UC} problem seeks minimum cost strategies to determine the commitment statuses of thermal generation resources (\textit{TGR}s) based on expected load, equipment limitations, and operational policies. The equipment limitations of \textit{TGR}s render \textit{UC} a time-coupled problem and require that the \textit{UC} decisions be taken typically one-hour to one-day ahead of dispatch based on the data/information available at the time of decision.\par
A key problem associated with \textit{UC} is the inherent uncertainty in load forecasts that \textit{GO}s need to rely on for decision making. As a ramification of the deepening penetration of \textit{RR}s, this uncertainty is compounded by the highly time-varying, intermittent, and uncertain nature of \textit{RR} power outputs. As such, effective short-term planning of power systems calls for \textit{UC} methodologies that undertake an explicit assessment of the uncertainty associated with net load.\par
Most studies in the literature leverage stochastic optimization (\textit{SO}) or robust optimization (\textit{RO}) techniques so as to address the uncertainty in net load in the \textit{UC} problem. Under the \textit{SO} paradigm, uncertain parameters are analytically characterized by a probability distribution that is specified \textit{a priori}. In fact, the optimal decisions are taken so as to minimize the expected cost based on the specified probability distribution. \par
\begin{table*}
\begin{minipage}[c][7.7cm][t]{\linewidth}
\begin{fminipage}{\textwidth}[1ex]
\footnotesize
\renewcommand{\arraystretch}{1.2}
\begin{tabular}{l l  l l }
\multicolumn{2}{l}{\textbf{Nomenclature}}  & & \\
$\mathscr{H}$/$h$  & set/index of simulation time periods	& $p_{g}$ & power generation of \textit{TGR} $g$ in hour $h$\\
$\mathscr{G}$/$g$ & set/index of thermal generation	& $p_{c}[h]$& curtailed load in hour $h$\\
& resources (\textit{TGR}s)	& $p_{s}[h]$& spilled \textit{RR} generation in hour $h$\\
$p_{g} ^{m}$/$p_{g} ^{M}$ & minimum/maximum power output	& $\boldsymbol{x}$ & vector of first-stage variables \\
 & of \textit{TGR} $g$ 	&&comprising $u_g[h]$ and $v_g[h]$\\
$T^{\uparrow}_{g} $/$T^{\downarrow}_{g}$ & minimum uptime/downtime of \textit{TGR} $g$	& $\boldsymbol{y}$ & vector of second-stage variables\\
$\Delta^{\uparrow}_{g} $/$\Delta^{\downarrow}_{g}$ & maximum ramp up/down rate of	&&comprising $p_g[h]$, $p_c[h]$, and $p_s[h]$\\
&  \textit{TGR} $g$	&$(\Omega, \mathcal{F}, \mathcal{P})$ & probability space\\
$\overline{\Delta}^{\uparrow}_{g} $/$\overline{\Delta}^{\downarrow}_{g}$ & start-up ramp up/shut-down ramp down 	&$\mathcal{P}_{o}$ & nominal probability distribution\\
& rate limit of \textit{TGR} $g$	&$\pi^{\omega}_{\mathcal{P}_{o}}$ & probability assigned to scenario $\omega$\\
$\lambda_{g}^{p}$/$\lambda_{g}^{u}$ & linear/fixed generation cost of \textit{TGR} $g$	&& by the distribution $\mathcal{P}$\\
$\lambda_{g}^{v}$ & start-up cost of \textit{TGR} $g$	&$\mathscr{P}$ & ambiguity set of probability distributions\\
$\lambda^{c}$ & cost of load curtailment	&$\boldsymbol{\tilde{\xi}}$ & random vector associated with net load  \\
$u_{g}[h]$/$v_{g}[h]$ & binary comitment status/start-up 	&$\boldsymbol{\xi^{\omega}}$ & realization $\omega$ of $\boldsymbol{\tilde{\xi}}$\\
&variable of  \textit{TGR} $g$ in hour $h$ 	&$\rho$  & divergence tolerance\\
\end{tabular}
\end{fminipage}
\end{minipage}
\end{table*}
While \textit{SO} suffers from a high computational burden, it delivers a lower total cost compared with its deterministic counterpart. This notwithstanding, the feasibility of \textit{SO} in real-life applications is hampered on two accounts. First, since \textit{GO}s have access to collected measurements but clearly not to their underlying distribution, they must commit to an assumed distribution in an \textit{SO} approach. Nevertheless, \textit{SO} renders a poor performance when the assumed distribution is wrong, and thus so does a pure \textit{SO} approach to the \textit{UC} problem. Second, when an \textit{SO} problem is solved based on a distribution tailored to a particular dataset and its out-of-sample performance is assessed on a different dataset, the obtained out-of-sample performance may be markedly low even when the two datasets are generated from the same distribution---a phenomenon referred to as the \textit{optimizer’s curse} \cite{druc:119}.\par
In contrast to \textit{SO} that considers the uncertain parameters to be of probabilistic nature, \textit{RO} models the uncertain parameters to vary within a predefined deterministic uncertainty set. By adopting a preventive viewpoint, \textit{RO} seeks to minimize the total cost based on the worst realization in the uncertainty set. While \textit{RO} provides an ex-ante protection against the worst-case realization, it has been subject to criticism for being liable to yield overly conservative decisions \cite{druc:109}.\par
In recent years, distributionally robust optimization (\textit{DRO}), also referred to as data-driven optimization, has gained prominence as a paradigm that capitalizes on the key strengths of both \textit{SO} and \textit{RO}. \textit{DRO} does not require that the underlying probability distribution be defined a priori, but rather relies on an ambiguity set of probability distributions. The hallmark of \textit{DRO} is the minimization of the expected cost induced by the worst-case distribution in the ambiguity set.\par
The notion of ambiguity set is the basic ingredient of the \textit{DRO} paradigm. The ambiguity set is constructed to comprise all probability distributions that are sufficiently close, as per the ambiguity set construction approach, to the nominal probability distribution, where the nominal probability distribution may be formed based on historical data, observations, or simulations. Approaches to ambiguity set construction may employ a measure to quantify the similarity between probability distributions, such as the Kullback-Leibler (\textit{KL}) divergence \cite{druc:109} and the Wasserstein metric \cite{druc:112}, or may use moment information \cite{druc:113}. To bound the distance/divergence between the nominal probability distribution and the probability distributions that are incorporated in the ambiguity set, these approaches typically rely on a tolerance parameter, such as the divergence tolerance in the \textit{KL} divergence or the radius of the Wasserstein ball in the Wasserstein metric. \par
Among these approaches, the \textit{KL} divergence, which shall be the approach employed in this paper, enjoys an asymptotic level of confidence that the true underlying probability distribution is contained within the constructed ambiguity set. This level of confidence is favorably impacted by rising divergence tolerance values, since they permit the incorporation of probability distributions that are farther from the nominal distribution, as well as by growing number of data points utilized in the formation of the nominal distribution.\par
In effectuating judicious \textit{UC} decisions, \textit{GO}s may be greatly aided by a distributionally robust unit commitment (\textit{DRUC}) methodology that expressly recognizes the uncertainty associated with net load---all the while hedging the \textit{UC} decisions against the inaptitude of the nominal probability distribution. The pursuance of such a data-driven methodology is further reinforced by the capability to collect vast amounts of data imparted by the efforts to modernize the legacy grid. The purpose of the present paper is the development of a methodology for \textit{DRUC} under net load uncertainty.\par
\subsection{Related Work}
In the literature, there has been an increasing preponderance of \textit{DRO} applications for the \textit{UC} problem over the last years. A \textit{DRO} framework that assesses the uncertainty in wind forecasting errors is worked out in \cite{druc:112}, where the devised framework capitalizes on the Wasserstein metric in constructing the ambiguity set. The authors of \cite{druc:113} propose a \textit{DRO} model for \textit{UC}, where they utilize moment information to build an ambiguity set. Nevertheless, the constructed ambiguity set relies on samples generated from a prespecified normal distribution, not actual measurements. \par
The \textit{DRO} model introduced in \cite{druc:114}  takes into account the uncertainty associated with \textit{RR} generation and load forecast errors in \textit{UC}, yet \cite{druc:114} leverages samples drawn from parametric probability distributions in constructing the ambiguity sets. The uncertainty in wind forecast errors has been studied in \cite{druc:109} under a \textit{DRO} framework involving an ambiguity set constructed using the \textit{KL} divergence. In the same vein as \cite{druc:113, druc:114}, the proposed model in \cite{druc:109} has been demonstrated using Monte Carlo samples taken from a normal distribution. 
\subsection{Contributions and Structure of the Paper}
The general contributions of this paper are as follows:
\begin{enumerate}
\item We develop a novel data-driven methodology for \textit{DRUC} under net load uncertainty. The proposed methodology may lend itself as a valuable tool for \textit{GO}s in the assessment of \textit{UC} decisions for systems with deep penetration of \textit{RR}s.
\item We lay out an ambiguity set construction approach that effectively exploits historical net load data. The presented approach is an extension of our previously reported work \cite{druc:116}, and it makes use of the k-means clustering algorithm in forming the nominal probability distribution and the \textit{KL} divergence in constructing the ambiguity set. We investigate the influence of different distance measures leveraged by the k-means clustering algorithm on the optimal solutions. 
\item By drawing on the performed empirical studies, we provide valuable insights into the sensitivity of the optimal solutions to dataset size and divergence tolerance. 
\end{enumerate}
The remainder of the paper consists of four sections. We map out the mathematical formulations of the proposed \textit{DRUC} methodology in Section \ref{2} and expound on the construction of the ambiguity set. In Section \ref{3}, we set forth a two-level decomposition algorithm for the solution of the devised \textit{DRUC} problem. We carry out representative studies in Section \ref{4} to illustrate the application of the proposed methodology and discuss the results. We present our concluding remarks in Section \ref{5}.
\section{Mathematical Underpinnings of the DRUC Methodology}\label{2}
In this section, we set up the mathematical formulation of the proposed Kullback-Leibler divergence-based distributionally robust unit commitment ($\mathrm{KL-DRUC}$) problem and present our ambiguity set construction approach.\par
The concept of time will play a predominant role in what is to follow, so it seems advisable to start out by describing our time notation. We discretize the time-axis and---commensurate with the typical time granularity and planning horizon of \textit{UC} studies---adopt 1 hour as the smallest indecomposable unit of time and 24 hours as the scheduling horizon. We define the study period by the set $\mathscr{H} \coloneqq \{h \colon h=1,...,24\}$.
\subsection{Problem Formulation}\label{2a}
In our $\mathrm{KL-DRUC}$ problem formulation, we expressly take into account the uncertainty associated with net load. We represent the uncertain net load over the study period by the random vector $\boldsymbol{\tilde{\xi}}$ defined on the probability space $(\Omega, \mathcal{F}, \mathcal{P})$, where $\Omega$ is a sample space,
$\mathcal{F}$ is a set of subsets of $\Omega$ that is a $\sigma-$algebra, and $\mathcal{P}$ is a probability distribution on $\mathcal{F}$. We represent by $\Xi$ the support of the probability distribution $\mathcal{P}$ and write the relation $\boldsymbol{\tilde{\xi}} \in \Xi \subset \mathbb{R}^{24}$. $\mathcal{P}$ is presumed to have a finite support with $\mathcal{S}$ realizations that we equivalently refer to as scenarios, i.e., $\vert\Xi\vert = \mathcal{S} < \infty$. We elaborate on the construction of $\Xi$ in Section \ref{2b}. \par
We denote by ${\xi^{\omega}}[h] $ the row $h$ of $\boldsymbol{\xi^{\omega}}$, which represents the net load in hour $h \in \mathscr{H}$ for scenario $\omega$. We denote by $\pi^{\omega}_{\mathcal{P}}$ the probability assigned to the scenario $\omega$ by the probability distribution $\mathcal{P}$. \par
The proposed $\mathrm{KL-DRUC}$ formulation involves two stages that reflect the sequence in which the \textit{UC} and economic dispatch decisions are rolled out. The first-stage problem is expressed as:
\begin{IEEEeqnarray}{ll}
\underset{u_g[h], v_g[h]}{\text{minimize}} &  \bigg\{\sum_{h \in \mathscr{H}}  \sum_{g \in \mathscr{G}} \Big[\lambda^{v}_{g} v_{g}[h] + \lambda^{u}_{g} u_{g}[h] \Big]  \nonumber\\
&\hspace{30pt}+ \,\underset{\mathcal{P} \in \mathscr{P}}{\text{maximize}} \;\mathbb{E}_{\mathcal{P}}\big[\mathcal{Q}(\boldsymbol{x},\boldsymbol{\tilde{\xi}})\big] \bigg\}, \label{obj}\\
\text{subject to} & \nonumber\\
&\hspace{-21pt}v_{g}[h] \geq u_{g}[h]-u_{g}[h-1] \;\; \forall g \in \mathscr{G}, \forall h \in \mathscr{H},\label{st}\\
&\hspace{-21pt}u_{g}[h]  -  u_{g}[h-1]  \leq  u_{g}[\tau]\;\; \forall \tau \in \mathbb{N}\: \text{such that} \nonumber\\
&\hspace{-21pt}\quad \quad  h \leq \tau \leq \min\{h-1+T^{\uparrow}_{g},24\}\;\forall g \in \mathscr{G}, \label{upt}  \\
&\hspace{-21pt}u_{g}[h-1]  -  u_{g}[h]  \leq  1-u_{g}[\tau]\;\; \forall \tau \in \mathbb{N}\: \text{such that}  \nonumber\\
&\hspace{-21pt}\quad \quad  h \leq \tau \leq \min\{h-1+T^{\downarrow}_{g},24\}\;\;\forall g \in \mathscr{G}, \label{dwt}  \\
&\hspace{-21pt} u_{g}[h], v_{g}[h] \in \{0,1\} \;\; \forall g \in \mathscr{G}, \forall h \in \mathscr{H}. \label{bin}
\end{IEEEeqnarray}
The first-stage problem evaluates the binary commitment statuses ($u_g[h]$) and start-up variables ($v_g[h]$) of the \textit{TGR}s over the study period with the objective \eqref{obj} to minimize the commitment and start-up costs plus the worst-case expected dispatch and load curtailment costs. The first-stage decisions are taken before the uncertain net load values get revealed subject to the minimum uptime \eqref{upt} and downtime \eqref{dwt} constraints of the \textit{TGR}s. We represent all first-stage decision variables by the vector $\boldsymbol{x}$, which comprises $u_g[h]$ and $v_g[h]$.\par
The inner maximization problem in \eqref{obj} is a key pillar in the $\mathrm{KL-DRUC}$ formulation, in that it ensures that the expected value of the variable generation and load curtailment costs be evaluated based on the worst-case probability distribution in the ambiguity set $\mathscr{P}$. As shall be spelled out in Section \ref{2b}, by virtue of this formulation, the \textit{DRUC} methodology has the ability to represent the probabilistic nature of uncertain net load, without being circumscribed by the characteristics of a single predefined probability distribution. \par
The function $\mathcal{Q}(\boldsymbol{x},\boldsymbol{\tilde{\xi}})$ in \eqref{obj} denotes the uncertain power generation and load curtailment costs. For a specific vector of first-stage decision variables $\boldsymbol{x}^{\dagger}$ and a realization $\boldsymbol{\xi^{\omega}}$, $\mathcal{Q}(\boldsymbol{x^{\dagger}},\boldsymbol{{\xi}^{\omega}})$ is computed by solving the following second-stage problem: 
\begin{IEEEeqnarray}{cl}
\mathcal{Q}(\boldsymbol{x^{\dagger}},\boldsymbol{{\xi}^{\omega}})\, \coloneqq   &\nonumber\\
 \underset{p_g[h], p_c[h], p_s[h]}{\text{minimize}} &  \sum_{h \in \mathscr{H}}\Big[ \sum_{g \in \mathscr{G}} {\lambda}^{p}_{g}p_{g}[h] + \lambda^{c}p_{c}[h]\Big], \label{obj2}\\
 \text{subject to} & \nonumber\\
&\hspace{-30pt}u_{g}^{\dagger}[h]p_{g}^{m} \leq p_{g}[h] \leq u_{g}^{\dagger}[h]p_{g}^{M} \nonumber\\ 
& \hspace{83pt} \forall g \in \mathscr{G}, \forall h \in \mathscr{H},\label{ol1} \\
& \hspace{-30pt}p_{g}[h] - p_{g}[h-1] \leq \Delta^{\uparrow}_{g}u_{g}^{\dagger}[h-1]  \nonumber\\
&\hspace{-10pt}+ \overline{\Delta}^{\uparrow}_{g}(1-u_{g}^{\dagger}[h-1])\quad\forall g \in \mathscr{G}, \forall h \in \mathscr{H},\label{ru}\\
& \hspace{-30pt}p_{g}[h-1] - p_{g}[h] \leq \Delta^{\downarrow}_{g}u_{g}^{\dagger}[h]  \nonumber\\
&\hspace{-10pt}+ \overline{\Delta}^{\uparrow}_{g}(1-u_{g}^{\dagger}[h])\hspace{28pt}\forall g \in \mathscr{G}, \forall h \in \mathscr{H},\label{rd}\\
&\hspace{-30pt}\sum_{g \in \mathscr{G}} p_{g}[h] + p_{c}[h] - p_{s}[h] = \xi^{\omega}[h] \quad \forall h \in \mathscr{H}, \label{pb}\\
&\hspace{-30pt}p_c[h],  p_s[h] \geq 0 \hspace{85pt} \forall h \in \mathscr{H}. \label{nn}
\end{IEEEeqnarray}
The second-stage problem \eqref{obj2}-\eqref{nn} seeks to determine the power dispatch of the \textit{TGR}s ($p_g[h]$), curtailed load ($p_c[h]$), and the spilled \textit{RR} generation ($p_s[h]$) over the study period. The objective of the second-stage problem \eqref{obj2} is to minimize the dispatch costs of \textit{TGR}s and the penalty cost incurred due to load curtailment, while taking into account the \textit{TGR} output limits \eqref{ol1}, ramping constraints \eqref{ru}-\eqref{rd}, and power balance constraint \eqref{pb}.
\subsection{Ambiguity Set Construction Approach}\label{2b}
We devote this subsection to the description of our approach to constructing $\mathscr{P}$. We start out by the construction of the nominal probability distribution $\mathcal{P}_{o}$. Let $\mathcal{N}$ denote the number of days for which historical net load data are initially considered. We leverage the k-means clustering algorithm to partition the $\mathcal{N}$ number of time-series data points to $\mathcal{S}$ clusters so as to assign each data point to the cluster with the nearest cluster centroid. To this end, distance measures such as the Euclidean distance (\textit{ED}), dynamic time warping (\textit{DTW}) \cite{druc:122}, or soft dynamic time warping (\textit{SDTW}) \cite{druc:106} may be harnessed by the k-means clustering algorithm in measuring the similarity between time-series data points, which may influence the constructed clusters as will be borne out in Section \ref{4}. \par
We utilize each of the constructed $\mathcal{S}$ clusters to form each of the $\mathcal{S}$ scenarios of the nominal probability distribution. For each cluster $\omega$, we draw on the cluster centroid to represent the realization $\boldsymbol{\xi^{\omega}}$ and construct the support $\Xi\coloneqq\{\boldsymbol{\xi^{\omega}}\colon\omega=1,\ldots,\mathcal{S}\}$. Let $\mathcal{N}^{\omega}$ denote the number of data points assigned to cluster $\omega$ by the k-means clustering algorithm. For the nominal probability distribution $\mathcal{P}_{o}$, we assign the probability for the scenario $\omega$ as $\pi^{\omega}_{\mathcal{P}_{o}}=\frac{\mathcal{N}^{\omega}}{\mathcal{N}},\,\omega=1,\ldots,\mathcal{S}$. \par
We next leverage the \textit{KL} divergence to construct an ambiguity set of probability distributions $\mathscr{P}$ around the nominal probability distribution $\mathcal{P}_{o}$. The ambiguity set formulation based on the \textit{KL} divergence \cite{druc:108} is expressed as:
\begin{IEEEeqnarray}{ll}
\vspace{0.17cm}
\mathscr{P}\coloneqq\{\mathcal{P} : & \sum_{\omega=1}^{S}\pi^{\omega}_{\mathcal{P}}\log\bigg(\frac{\pi^{\omega}_{\mathcal{P}}}{\pi^{\omega}_{\mathcal{P}_{o}}}\bigg)\leq\rho, \label{dc}\\
& \sum_{\omega=1}^{S}\pi^{\omega}_{\mathcal{P}}=1, \label{sc}\\
& \pi^{\omega}_{\mathcal{P}}\geq0\quad\forall \omega \in \Omega\}. \label{pc}
\vspace{0.17cm}
\end{IEEEeqnarray}
The size of the ambiguity set is governed by the divergence tolerance $\rho$ through the relation \eqref{dc}. Setting $\rho=0$ yields an ambiguity set $\mathscr{P}$ that includes solely the nominal distribution and thus strips $\mathscr{P}$ of its ability to immunize against adopting an unsuitable probability distribution. At the other extreme, as $\rho \rightarrow \infty$, $\mathscr{P}$ admits all probability distributions. Our empirical studies in Section \ref{4} will lay bare the influence of $\rho$ on optimal solutions. While increasing the value of $\rho$ drives up the worst-case expected costs in \eqref{obj} and may lead to overly conservative decisions, it also results in larger ambiguity sets that will typically contain the true underlying probability distribution with a greater level of confidence. \par
In \cite{druc:115}, Ben-Tal et al. work out an asymptotic value for the divergence tolerance $\rho$ that delivers an approximate confidence level of $1-\eta$ for having included the true underlying probability distribution in the constructed ambiguity set. This asymptotic value of $\rho$ is expressed as: 
\begin{IEEEeqnarray}{l}
\vspace{0.17cm}
\rho=\frac{1}{2\mathcal{N}}\chi^2_{\mathcal{S}-1, 1-\eta}, \label{asymp}
\vspace{0.17cm}
\end{IEEEeqnarray}
where $\chi^2_{\mathcal{S}-1, 1-\eta}$ denotes the $1-\eta$ quantile of the $\chi^2_{\mathcal{S}-1}$ distribution with $\mathcal{S}-1$ degrees of freedom. The expression in \eqref{asymp} is based on asymptotics and so is only approximately valid. Based on \eqref{asymp}, for a fixed number of samples $\mathcal{N}$, increasing divergence tolerance values bring about a greater level of confidence that the true underlying probability distribution is included within the ambiguity set, which conforms with our earlier remarks. Further, with the utilization of greater number of samples $\mathcal{N}$, a fixed level of confidence tolerance $1-\eta$ can be attained by a smaller divergence tolerance $\rho$ and hence a smaller ambiguity set, which analytically affirms the value of collecting additional data. This analytical relation will be corroborated by our studies in Section \ref{4}. \par
We succinctly state the $\mathrm{KL-DRUC}$ formulation as:
\begin{IEEEeqnarray}{ll}
 \hspace{-15pt}\mathrm{KL-DRUC:}   &\nonumber\\
 \hspace{-15pt}\underset{\boldsymbol{x}}{\text{minimize}} & \hspace{7pt} \boldsymbol{c}\cdot\boldsymbol{x} +\underset{\mathcal{P} \in \mathscr{P}}{\text{maximize}} \;\sum_{\omega \in \Omega}\pi^{\omega}_{\mathcal{P}}\mathcal{Q}(\boldsymbol{x},\boldsymbol{{\xi}^{\omega}}),  \label{obj3}\\
 \hspace{-15pt}\text{subject to} & \hspace{7pt}x \in \mathscr{X} , \label{xc}\\
& \hspace{7pt} \eqref{dc}-\eqref{pc}, \nonumber
\end{IEEEeqnarray}
where $\mathscr{X}$ represents the feasibility region of $\boldsymbol{x}$ defined by the constraints \eqref{st}-\eqref{bin}.
\section{Solution Method}\label{3}
We devote this section to the development of a method based on Benders' decomposition for the efficient solution of the $\mathrm{KL-DRUC}$ problem. Since the $\mathrm{KL-DRUC}$ problem has an intractable min-max-min structure, we leverage the duality theory so as to recast the $\mathrm{KL-DRUC}$ problem into a tractable form. By assigning the dual variables $\zeta$ and $\mu$ to the constraints \eqref{dc} and \eqref{sc}, respectively, we take the dual of the inner maximization problem in $\mathrm{KL-DRUC}$ and obtain the following convex mixed-integer nonlinear reformulated $\mathrm{KL-DRUC}$ ($\mathrm{RKL-DRUC}$) problem:
\begin{IEEEeqnarray}{ll}
 \hspace{-15pt}\mathrm{RKL-DRUC}:  &\hspace{-3pt}\nonumber\\
 \hspace{-15pt}\underset{\boldsymbol{x}, \mu, \zeta}{\text{minimize}} & \hspace{-3pt} \boldsymbol{c}\cdot\boldsymbol{x} + \mu + \rho \zeta + \zeta\sum_{\omega =1}^{\mathcal{S}} \pi^{\omega}_{\mathcal{P}_{o}}e^{\overline{\mathcal{K}}^{\omega}(\boldsymbol{x}, \mu, \zeta)-1},  \label{obj4}\\
 \hspace{-15pt}\text{subject to} & \hspace{-3pt}x \in \mathscr{X} , \label{xc}\\
 \hspace{0pt}& \hspace{-3pt}\zeta \geq 0, \label{lc}
\end{IEEEeqnarray}
where $\overline{\mathcal{K}}^{\omega}(\boldsymbol{x}, \mu, \zeta)=\frac{\mathcal{Q}(\boldsymbol{x},\boldsymbol{\xi^{\omega}})-\mu}{\zeta}$ \cite{druc:108}. For notational brevity, we define the following functions: 
\begin{IEEEeqnarray}{ll}
{\overline{\mathcal{R}}^{\omega}}(\boldsymbol{x},\zeta, \mu)  & \coloneqq  \zeta e^{\overline{\mathcal{K}}^{\omega}(\boldsymbol{x}, \mu, \zeta)-1},\label{rbarfnc}\\
{\mathcal{R}}(\boldsymbol{x},\zeta, \mu) & \coloneqq \sum_{\omega =1}^{\mathcal{S}} \pi^{\omega}_{\mathcal{P}_{o}}{\overline{\mathcal{R}}^{\omega}}(\boldsymbol{x},\zeta, \mu).\label{rfnc}
\end{IEEEeqnarray}
Benders' decomposition involves the solution of an optimization problem with complicating variables in a distributed manner at the expense of iterations. By capitalizing on Benders' decomposition, we solve a lower-bounding master problem ($\mathrm{MP}$) and an upper-bounding subproblem ($\mathrm{SP}$) in an iterative way, rather than solving the monolithic $\mathrm{RKL-DRUC}$ problem for all decisions variables simultaneously. As per \cite{druc:108}, we consider $\boldsymbol{x}$, $\zeta$, and $\mu$ to be the complicating variables and, for a Benders' iteration $\nu$, set up the following $\mathrm{MP}$:\\
\begin{IEEEeqnarray}{ll}
 \hspace{0pt}  \mathrm{MP:}   &\hspace{6pt}\nonumber\\
 \hspace{0pt} \underset{\boldsymbol{x_{(\nu)}}, \mu_{(\nu)}, \zeta_{(\nu)}}{\text{minimize}} &\hspace{8pt}  \boldsymbol{c}\cdot\boldsymbol{x_{(\nu)}} + \mu_{(\nu)} + \rho \zeta_{(\nu)} + \theta_{(\nu)} \label{obj5}\\
\hspace{0pt}  \text{subject to} &\hspace{6pt} x_{(\nu)} \in \mathscr{X} , \label{xc}\\
 \hspace{0pt} &\hspace{6pt} \zeta_{(\nu)} \geq 0, \label{lc}\\
 \hspace{0pt}  & \hspace{6pt}\theta_{(\nu)} \geq \boldsymbol{{\alpha}_{(j)}}\cdot  (\boldsymbol{x_{(\nu)}}-\boldsymbol{x_{(j)}}) + {\beta}_{(j)} (\mu_{(\nu)}- \mu_{(j)})\nonumber\\
\hspace{35pt} & \hspace{35pt}+  {\gamma}_{(j)} (\zeta_{(\nu)}-\zeta_{(j)}) + {\mathcal{R}}(\boldsymbol{x_{(j)}}, \zeta_{(j)}, \mu_{(j)}), \nonumber\\
 \hspace{35pt}&\hspace{35pt}\,j =1,\ldots,\nu-1. \label{optc}
\end{IEEEeqnarray}\par
The $\mathrm{MP}$ is a relaxed version of the $\mathrm{RKL-DRUC}$ problem, in that the objective function of the $\mathrm{MP}$ \eqref{obj5} approximates from below the objective function of the $\mathrm{RKL-DRUC}$ \eqref{obj4}. This approximation is imparted by the Benders' optimality cuts expressed in \eqref{optc}, which serve to approximate from below the function $\mathcal{R}(\boldsymbol{x},\zeta, \mu)$. As the number of iterations increases, more and more Benders' optimality cuts are generated and evaluated in \eqref{optc}, thereby gradually rendering the $\mathrm{MP}$ less relaxed.\par
At each iteration $\nu$, the optimal $\mathrm{MP}$ solution $(\boldsymbol{x_{(\nu)}}, {\mu_{(\nu)}}, {\zeta_{(\nu)}})$ is fixed as $\boldsymbol{x_{f}}\gets\boldsymbol{x_{(\nu)}},$
${\mu_{f}}\gets {\mu_{(\nu)}}$, and ${\zeta_{f}}\gets{\zeta_{(\nu)}}$ and utilized in solving the $\mathrm{SP}$, which is described next.\par
In line with \cite{druc:116, druc:108}, we adopt $\mathcal{Q}(\cdot)$ presented in \eqref{obj2}-\eqref{nn} as the $\mathrm{SP}$ and harness the chain rule along with the optimal $\mathrm{SP}$ solution in evaluating the optimality cuts. At a given iteration $\nu$, we express the $\mathrm{SP}$ for scenario $\omega$ as: 
\begin{IEEEeqnarray}{lll}
 \hspace{-15pt}\mathrm{SP:}   &\hspace{7pt} &\nonumber\\
 \hspace{-15pt}\underset{\boldsymbol{y^{\omega}_{(\nu)}}, \boldsymbol{\hat{x}^{\omega}_{(\nu)}}}{\text{minimize}} &  \hspace{7pt} \sum_{h \in \mathscr{H}}  \sum_{g \in \mathscr{G}} \lambda_{g}^{p}p_{g}[h] + \lambda^{c}p_{c}[h], \label{obj22}\\
 \hspace{-15pt} \text{subject to} & \hspace{7pt} \eqref{ol1}-\eqref{nn},\nonumber\\
 & \hspace{7pt} \boldsymbol{\hat{x}^{\omega}_{(\nu)}}= \boldsymbol{x_{f}}& \hspace{-50pt}\xleftrightarrow{}\hspace{20pt} \boldsymbol{\overline{\varphi}^{\omega}_{(\nu)}}.\label{alc}
 \end{IEEEeqnarray}
The dual variable $ \boldsymbol{\overline{\varphi}^{\omega}_{(\nu)}}$ associated with the constraint \eqref{alc} denotes the negative of the sensitivity of \eqref{obj22} to $\boldsymbol{x_{f}}$. Note that the elements of $\boldsymbol{\hat{x}^{\omega}_{(\nu)}}$ are not constrained to be binary and so the $\mathrm{SP}$ for each scenario $\omega$ is a continuous problem. To compute the Benders' optimality cuts, we evaluate the following terms for each scenario $\omega$:
\begin{IEEEeqnarray}{ll}
\boldsymbol{\overline{\alpha}^{\omega}_{(\nu)}}&=\zeta_{f} e^{\overline{\mathcal{K}}^{\omega}(\boldsymbol{x_{f}}, \mu_{f}, \zeta_{f})-1}\, \boldsymbol{\overline{\varphi}^{\omega}_{(\nu)}}, \label{alphaeq}\\
\overline{\beta}^{\omega}_{(\nu)}&=\frac{\partial {\overline{\mathcal{R}}^{\omega}}(\boldsymbol{x_{f}},\zeta_f, \mu_f)}{\partial \zeta_f}\nonumber\\
&=(1-\overline{\mathcal{K}}^{\omega}(\boldsymbol{x_{f}}, \mu_{f}, \zeta_{f}))e^{\overline{\mathcal{K}}^{\omega}(\boldsymbol{x_{f}}, \mu_{f}, \zeta_{f})-1},\label{betaeq}	\\
\overline{\gamma}^{\omega}_{(\nu)}&=\frac{\partial {\overline{\mathcal{R}}^{\omega}}(\boldsymbol{x_{f}}, \zeta_f, \mu_f)}{\partial \mu_f}=-e^{\overline{\mathcal{K}}^{\omega}(\boldsymbol{x_{f}}, \mu_{f}, \zeta_{f})-1},\label{gammaeq}
\end{IEEEeqnarray}
and compute the terms $\boldsymbol{{\alpha}_{(\nu)}}$, ${{\beta}_{(\nu)}}$, ${{\gamma}_{(\nu)}}$ in the Benders' optimality cuts in \eqref{optc} as follows: $\boldsymbol{{\alpha}_{(\nu)}}=\sum_{\omega =1}^{\mathcal{S}} \pi_{\mathcal{P}_{o}}^{\omega} \boldsymbol{\overline{\alpha}^{\omega}_{(\nu)}}$, ${{\beta}_{(\nu)}}=\sum_{\omega =1}^{\mathcal{S}} \pi_{\mathcal{P}_{o}}^{\omega} {\overline{\beta}^{\omega}_{(\nu)}}$, ${{\gamma}_{(\nu)}}=\sum_{\omega=1}^{\mathcal{S}} \pi_{\mathcal{P}_{o}}^{\omega} {\overline{\gamma}^{\omega}_{(\nu)}}$ \cite{druc:116}. We succinctly express the proposed method in Algorithm \ref{algo}.\par
 \vspace{-0.2cm}
\begin{algorithm}
 \vspace{-0cm}
\caption{Decomposition algorithm for $\mathrm{RKL-DRUC}$}\label{algo}
  \begin{algorithmic}[1]
    \STATE \text{Initialize} $\boldsymbol{x}\gets\boldsymbol{0}$. 
    \STATE \text{Solve} $\mathrm{SP}$. \text{Set} $\mathcal{Q}^{M}\gets \max(\{\mathcal{Q}(\boldsymbol{x}, \boldsymbol{{\xi}^{\omega}})\colon \omega \in \Omega\})$
    \STATE \text{Initialize} $\mathrm{UB}\gets\infty$, $\mathrm{LB}\gets-\infty$, $\nu\gets1$, $\theta_{(1)}\gets0$, $\zeta_{(1)}\gets0$, $\boldsymbol{x_{(1)}}\gets\boldsymbol{0}$.
    \WHILE {$\mathrm{UB}-\mathrm{LB}\geq\mathrm{TOL}$}
    	\STATE \text{Solve} $\mathrm{MP}$. \text{Determine} $\boldsymbol{x_{(\nu)}}, \zeta_{(\nu)}, \mu_{(\nu)}$, and $\theta_{(\nu)}$ so that  $\frac{\mathcal{Q}^{M}-\mu_{(\nu)}}{\zeta_{(\nu)}}\leq\mathcal{K}^{M}$. $\mathrm{LB}\gets\theta_{(\nu)}$.
	\STATE \text{Solve} $\mathrm{SP}$. \text{Determine} $\boldsymbol{{\alpha}_{(\nu)}}$, ${\beta}_{(\nu)}$, ${\gamma}_{(\nu)}$, and ${\mathcal{R}}(\boldsymbol{x_{(\nu)}}, \zeta_{(\nu)}, \mu_{(\nu)})$. $\mathrm{UB}\gets {\mathcal{R}}(\boldsymbol{x_{(\nu)}}, \zeta_{(\nu)}, \mu_{(\nu)})$.  $\nu\gets \nu+1$. 
    \ENDWHILE
  \end{algorithmic}
\end{algorithm}
\section{Case Study and Results}\label{4}
In this section, we carry out representative studies to demonstrate the effectiveness of our proposed methodology. Our studies involve three \textit{TGR}s whose characteristics are based on \cite{druc:118}. The source code of the case studies are available at \cite{druc:111}. The aggregate peak capacity of the \textit{TGR}s is 1083 \textit{MW}. To construct our net load dataset, we use the net load values recorded in the California Independent System Operator grid from July 1, 2018 to October 31, 2020 \cite{druc:117} and scale the recorded net load values so that the maximum value in the net load dataset is scaled down to 1083 \textit{MW}. \par
In the foregoing sections, we laid out that the divergence tolerance $\rho$, the distance measure leveraged by the k-means clustering algorithm, and the dataset size may bear directly on the optimal solutions. Since it is impractical to traverse the entire set of possible values for these parameters, in what follows, we shall conduct representative studies to probe the influence of each parameter.\par
We initially consider the net load values reported for the first twelve months, that is, July 1, 2018-June 30, 2019, and leverage the k-means clustering algorithm separately with each of the \textit{ED}, \textit{DTW}, and \textit{SDTW} distance measures. To pick the number of clusters, we assess the percentage of variance captured under different values of $\mathcal{S}$ for each distance measure so as to identify the point of diminishing returns, i.e., the so-called \textit{elbow}. While the elbow for each distance measure is observed at a slightly different value of $\mathcal{S}$, we aim at the identification of a single cluster number that can be consistently used for all distance measures in order to be able to study the impact of the employed distance measure on a uniform basis. For this purpose, we pick the smallest $\mathcal{S}$ that can capture a percentage of variance greater than or equal to that of the elbow for all three distance measures, \textit{viz.}: $\mathcal{S}=12$. The constructed three sets of clusters are harnessed in constructing three nominal probability distributions and their associated supports.\par
We solve the $\mathrm{RKL-DRUC}$ problem by capitalizing on the solution method set forth in Section \ref{3}. We perform our implementations in \textit{Pyomo} using Gurobi 9.0.2 as the solver with the optimality tolerance gap $\mathrm{TOL}=10^{-4}$ on a 2.6 \textit{GHz} Intel Core i7 \textit{CPU} with 16 \textit{GB} of \textit{RAM}. Our discussion in Section \ref{2} drove home the pivotal role the divergence parameter $\rho$ plays in the degree of conservatism of ambiguity sets. To gain quantitative insights into the impact of $\rho$, we solve the $\mathrm{RKL-DRUC}$ problem by using each of the three constructed nominal probability distributions for each of the following values of $\rho=0,\,0.2,\,0.4,\,0.6,\,0.8,\,\text{and}\,1.0$.\par 
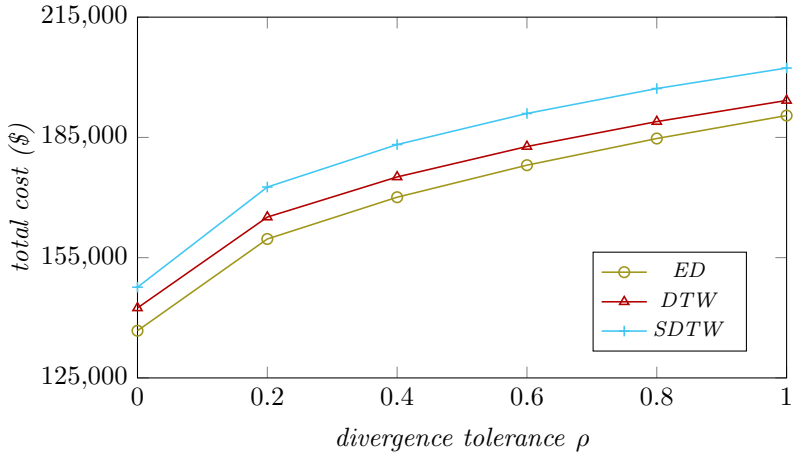
\begin{figure}[!h]
\vspace{0cm}
\centering
\begin{tikzpicture}
\pgfplotsset{
yticklabel style={
        /pgf/number format/fixed,
        /pgf/number format/precision=5
},
scaled y ticks=false}
\begin{axis}[
xlabel={\textit{divergence tolerance }$\rho$},
ylabel={\textit{total cost (\textit{\$})}},
width=0.64\textwidth,
height=0.40\textwidth,
xlabel style={at={(axis description cs:0.5, 0.0)},anchor=north},
ylabel style={at={(axis description cs:0.0,.5)},rotate=0, anchor=south},
xmin=0,
xmax=1.0,
xtick = {0,0.2,0.4,0.6, 0.8, 1.0},
ymin=125000,
ymax=215000,
ytick = {125000,155000, 185000, 215000},
legend style={font=\scriptsize, at={(0.82,0.35)}, legend columns=1, anchor=north},
label style={font=\small},
tick label style={font=\small}  ]
	\addplot+[sharp plot, white!10!olive, solid, semithick, mark=o] plot coordinates
		{ (0, 136771.4888383576) (0.2, 159659.7288582247) (0.4, 170076.0688727766) (0.6, 178068.92426174376) (0.8, 184718.9556282434) (1.0, 190451.98150559055)  };
	\addlegendentry{\textit{ED} \hspace{0.12cm}}
	\addplot+[sharp plot, black!30!red, solid, semithick, mark=triangle] plot coordinates
		{ (0, 142454.25684298886) (0.2, 165118.8558957899) (0.4, 175131.17902727475) (0.6, 182724.87707343895) (0.8, 188942.14181602825) (1.0, 194214.09467119823)  };
	\addlegendentry{\textit{DTW} \hspace{0.12cm}}
		\addplot+[sharp plot, white!40!cyan, solid, semithick, mark=+] plot coordinates
		{ (0, 147630.92939805114) (0.2, 172611.90759731026) (0.4, 183210.37295925722) (0.6, 190974.3466935816) (0.8, 197190.41380131827) (1.0, 202307.32576272113)  };
	\addlegendentry{\textit{SDTW} \hspace{0.12cm}}
	\end{axis}
\end{tikzpicture}
\vspace{0.0cm}
\caption{Optimal solutions under different distance measures as a function of divergence tolerance $\rho$}
\vspace{0.0cm}
\label{res1}
\end{figure}
Fig. \ref{res1} depicts the $\mathrm{RKL-DRUC}$ problem solutions for each distance measure as a function of the divergence tolerance $\rho$. We observe from Fig. \ref{res1} that, for all three distance measures, the total cost increases with larger values of $\rho$. This observation bears out the fact that larger $\rho$ values permit the incorporation of more adverse probability distributions into the ambiguity set, thereby increasing the optimal worst-case expected cost.\par 
We next turn to the investigation of the influence of dataset size on optimal total cost. To this end, we make use of the daily historical net load values reported in 15 time windows, where each time window has a different time span and thus encompasses a different number of time-series data points. All 15 time windows start from July 1, 2018 and extend to comprise the net load values for the following number of months: 1, 2, 4, 6, 8, 10, 12, 14, 16, 18, 20, 22, 24, 26, and 28. Using the net load values for each of these 15 time windows, we repeat the process for constructing the nominal probability distribution spelled out in Section \ref{2b} under each of the said 3 distance measures and ultimately obtain 45 nominal probability distributions and their associated supports.\par
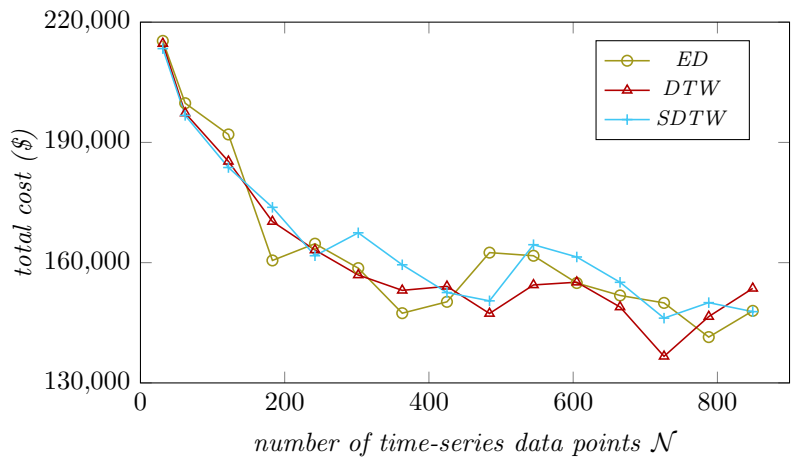
\begin{figure}[!h]
\vspace{0cm}
\centering
\begin{tikzpicture}
\pgfplotsset{
yticklabel style={
        /pgf/number format/fixed,
        /pgf/number format/precision=5
},
scaled y ticks=false}
\begin{axis}[
xlabel={\textit{number of time-series data points }$\mathcal{N}$},
ylabel={\textit{total cost (\textit{\$})}},
width=0.64\textwidth,
height=0.40\textwidth,
xlabel style={at={(axis description cs:0.5, 0.0)},anchor=north},
ylabel style={at={(axis description cs:0.0,.5)},rotate=0, anchor=south},
xmin=0,
xmax=900,
xtick = {0,200,400,600, 800},
ymin=130000,
ymax=220000,
ytick = {130000,160000, 190000, 220000},
legend style={font=\scriptsize, at={(0.82,0.95)}, legend columns=1, anchor=north},
label style={font=\small},
tick label style={font=\small}  ]
	\addplot+[sharp plot, white!10!olive, solid, semithick, mark=o] plot coordinates
		{ (31, 215323.8181793563) (62, 199781.88792846532) (122, 191984.3490963542) (183, 160550.77956499017) (242, 164732.61512652025) (302, 158629.81593159228)  (363, 147395.5965536313) (425, 150231.00236535523) (484, 162504.9253955754) (545, 161723.18681425243) (605, 154922.2652383383) (665, 151820.80369023595) (726, 149960.46657375334) (788, 141433.6533825267) (849, 147975.22397226543)};
	\addlegendentry{\textit{ED} \hspace{0.12cm}}
	\addplot+[sharp plot, black!30!red, solid, semithick, mark=triangle] plot coordinates
				{ (31, 214650.24690724773) (62, 197324.2473702508) (122, 185269.18224349705) (183, 170260.38384086173) (242, 163161.03873144492) (302, 156927.23823202794)  (363, 153092.93200074614) (425, 154089.66340378832) (484,  147316.30110467994) (545, 154442.0112323526) (605, 155124.40169356533) (665, 148938.85192655333) (726, 136630.50301633755) (788, 146548.61420400225) (849, 153582.24193353875)};
	\addlegendentry{\textit{DTW} \hspace{0.12cm}}
		\addplot+[sharp plot, white!40!cyan, solid, semithick, mark=+] plot coordinates
				{ (31, 213395.26983847065) (62, 196738.71710436422) (122, 183738.91380555456) (183, 173804.22046968905) (242, 161733.9396283884) (302, 167446.50749978525)  (363, 159470.52442052952) (425, 152549.4939670907) (484,  150477.11030025498) (545, 164480.19633756156) (605, 161418.97410521505) (665, 155102.60358840178) (726, 146161.35197446257) (788, 150025.66246000104) (849, 147778.89582505988)};
	\addlegendentry{\textit{SDTW} \hspace{0.12cm}}
	\end{axis}
\end{tikzpicture}
\vspace{0.0cm}
\caption{Optimal solutions under different distance measures as a function of the number of time-series data points $\mathcal{N}$}
\vspace{0.0cm}
\label{res2}
\end{figure}

For each of the 45 nominal probability distributions, we construct an ambiguity set whose divergence tolerance $\rho$ is picked so as to ensure a confidence level of $1-\eta=0.98$ as per \eqref{asymp}. By analogously employing the solution method illustrated in Section \ref{3}, we solve the $\mathrm{RKL-DRUC}$ problem for each of these 45 setups and present the results in Fig. \ref{res2}.\par
\enlargethispage{-0.3cm}
The plots in Fig. \ref{res2} demonstrate the tight coupling between the optimal cost and the number of data points utilized in constructing the nominal probability distribution. Independent of the distance measure, the optimal cost, by and large, decreases with growing number of data points. This decrease in optimal cost is in essential agreement with the asymptotic level of confidence brought forward in Section \ref{2}, which made clear that increasing values of $\mathcal{N}$ bring about smaller $\rho$ values for a fixed level of confidence $1-\eta$, thus shrinking the ambiguity set and driving down the worst-case expected cost.  
\section{Conclusion}\label{5}
In this paper, a distributionally robust unit commitment methodology using the Kullback-Leibler divergence has been worked out. The proposed methodology could aid grid operators in the evaluation of unit commitment decisions for systems with deep penetration of renewable resources. Grid operators may effectively exploit collected data jointly with our methodology so as to construct an ambiguity set of probability distributions, which affords the capability to represent the probabilistic nature of uncertain parameters without relying on a prespecified probability distribution. The cornerstone of the proposed methodology is the minimization of the total cost induced by the worst-case probability distribution in the ambiguity set. We carry out extensive numerical studies to demonstrate the effectiveness of the proposed methodology. Our sensitivity analyses shed light on the influence of divergence tolerance and dataset size on optimal total cost. 

\end{document}